\documentclass[11pt]{article}%
\usepackage{amsmath}
\usepackage{amsfonts}
\usepackage{amssymb}
\usepackage{graphicx}%
\providecommand{\U}[1]{\protect\rule{.1in}{.1in}}
\newtheorem{theorem}{Theorem}[section]

\newtheorem{remark}[theorem]{Remark}

\numberwithin{equation}{section}

\begin{document}

\title{Reconstruction techniques for complex potentials}
\author{Vladislav V. Kravchenko\\{\small Departamento de Matem\'{a}ticas, Cinvestav, Unidad Quer\'{e}taro, }\\{\small Libramiento Norponiente \#2000, Fracc. Real de Juriquilla,
Quer\'{e}taro, Qro., 76230 MEXICO.}\\{\small e-mail: vkravchenko@math.cinvestav.edu.mx}}
\maketitle

\begin{abstract}
An approach for solving a variety of inverse coefficient problems for the
Sturm-Liouville equation $-y^{\prime\prime}+q(x)y=\rho^{2}y$ with a complex
valued potential $q(x)$ is presented. It is based on Neumann series of Bessel
functions representations for solutions. With their aid the problem is reduced
to a system of linear algebraic equations for the coefficients of the
representations. The potential is recovered from an arithmetic combination of
the first two coefficients. Special cases of the considered problems include
the recovery of the potential from a Weyl function, inverse two-spectra
Sturm-Liouville problems, as well as the inverse scattering problem on a
finite interval. The approach leads to efficient numerical algorithms for
solving coefficient inverse problems. Numerical efficiency is illustrated by
several examples.

\end{abstract}

\section{Introduction}

Inverse coefficient problems for second order linear differential equations
arise in numerous applications in modern physics, chemistry, Earth sciences
and engineering (see, e.g., \cite{Chadan et al 1997}, \cite{Freiling and
Yurko}, \cite{Gladwell}, \cite{Gou Chen 2015}, \cite{Kabanikhin},
\cite{KrBook2020}, \cite{LevitanInverse}, \cite{Marchenko}, \cite{Ramm 2005},
\cite{Vatulyan 2019}, \cite{Yurko2007}). The theory of different classes of
such problems is especially well developed for the one-dimensional
Schr\"{o}dinger equation, which is also a Liouville normal form of a general
Sturm-Liouville equation
\begin{equation}
-y^{\prime\prime}(x)+q(x)y(x)=\rho^{2}y(x),\quad x\in(0,L).
\label{SL equation}%
\end{equation}
Here $\rho\in\mathbb{C}$ is a spectral parameter and the variable coefficient
$q\in%
\mathcal{L}%
_{2}(0,L)$ is in general a complex-valued function, called the potential. We
will consider equation (\ref{SL equation}) on a finite interval $(0,L)$, $L>0$.

A variety of inverse problems for (\ref{SL equation}) have been extensively
studied, such as the recovery of $q$ from a spectral function, or from two
spectra, or from a Weyl-Titchmarsh function (see, e.g., \cite{Chadan et al
1997}, \cite{Freiling and Yurko}, \cite{KrBook2020}, \cite{LevitanInverse},
\cite{Marchenko}, \cite{Yurko2007}). The theory of those and several other
inverse problems for (\ref{SL equation}), questions of existence and
uniqueness have been studied rather thoroughly. More recently the attention
has shifted to the question of obtaining efficient reconstructive methods, and
this paper continues the trend.

Meanwhile the theory of several kinds of inverse problems for
(\ref{SL equation}) is quite well developed even in the case when the
potential is complex valued (a non-selfadjoint case)
\cite{AvdoninBelishev1996}, \cite{Bondarenko2020}, \cite{Buterin2007},
\cite{ButerinKuznetsova2019}, \cite{Horvath}, \cite{Marletta},
\cite{Marchenko}, \cite{XuBondarenko2023}, to our best knowledge all the
developed numerical algorithms (see, e.g., \cite{Brown et al 2003}, \cite{Gao
et al 2013}, \cite{Gao et al 2014}, \cite{IgnatievYurko}, \cite{Kammanee
Bockman 2009}, \cite{Kr2019JIIP}, \cite{KrBook2020}, \cite{KrSpectrumCompl},
\cite{KKC2022Mathematics}, \cite{KST2020IP}, \cite{KT2021 IP1}, \cite{KT2021
IP2}, \cite{Neamaty et al 2017}, \cite{Neamaty et al 2019}, \cite{Rohrl},
\cite{Rundell Sacks}, \cite{Sacks}) require the potential to be real valued.
Among the difficulties arising in a non-selfadjoint situation we mention the
difficulty of recovering two spectra from the Weyl function or of applying the
Gelfand-Levitan integral equation approach. In particular, the requirement of
the real valuedness of the potential was essential in the recent publications
\cite{AKK2023tree}, \cite{AvdoninKravchenko2022}, \cite{KrSpectrumCompl},
\cite{KKC2022Mathematics}, because the approach developed there for solving a
variety of inverse problems, required computing two spectra for
(\ref{SL equation}), which is a challenging numerical problem when $q(x)$ is
complex. The purpose of this paper is to present a simple and efficient method
for the numerical solution of a variety of inverse problems for
(\ref{SL equation}), which recovers equally well both real and complex-valued potentials.

Consider the following model. Let $u(\rho,x)$ be a solution of
(\ref{SL equation}), satisfying some prescribed initial conditions at the
origin
\begin{equation}
u(\rho,0)=a(\rho),\quad u^{\prime}(\rho,0)=b(\rho), \label{init}%
\end{equation}
where $a(\rho)$ and $b(\rho)$ are some known functions, in general
complex-valued, which are not identical zeros. Denote $\ell(\rho)=u(\rho,L)$.
For any given $a(\rho)$ and $b(\rho)$, a unique solution $u(\rho,x)$ exists,
and so $\ell(\rho)$ is well defined. Let $\left\{  \rho_{k}\right\}
_{k=1}^{\infty}$ be an infinite set of complex numbers, and
\begin{equation}
a_{k}:=a(\rho_{k}),\quad b_{k}:=b(\rho_{k}),\text{\quad}\ell_{k}:=\ell
(\rho_{k}),\text{\quad}k=1,2,\ldots. \label{lk}%
\end{equation}
The main aim of the present work is to develop a method for the numerical
solution of the following inverse problem.

\textbf{Problem A.} Given the data
\begin{equation}
\left\{  \rho_{k},\,a_{k},\,b_{k},\,\ell_{k}\right\}  _{k=1}^{\infty},
\label{inpdata}%
\end{equation}
find $q(x)$.

As was pointed out in \cite{KKC2022Mathematics}, special cases of Problem A
include such well studied inverse problems as the recovery of the potential
from a Weyl function or from two spectra. Indeed,\ let $\Phi(\rho,x)$ denote a
Weyl solution of (\ref{SL equation}), which satisfies (\ref{SL equation})
together with the boundary conditions
\[
\Phi^{\prime}(\rho,0)=1\quad\text{and\quad}\Phi(\rho,L)=0.
\]
If $\rho^{2}$ is not a Neumann-Dirichlet eigenvalue of (\ref{SL equation}),
the solution $\Phi(\rho,x)$ exists and is unique. Let $M(\rho)$ denote the
Weyl function $M(\rho):=\Phi(\rho,0)$. Then one of the frequently studied
inverse Sturm-Liouville problems (see, e.g., \cite{Bondarenko1},
\cite{Bondarenko2}, \cite{Gesztesy}, \cite{Horvath}, \cite{Yurko2007})
consists in recovering $q(x)$ from $M(\rho)$ (given at a set of points). In
other words, we have
\begin{equation}
a_{k}=M(\rho_{k}),\ b_{k}=1\quad\text{and\quad}\ell_{k}=0,\ k=1,2,\ldots.
\label{dataWeyl}%
\end{equation}
From these data the potential $q(x)$ needs to be recovered.

A classical two-spectra inverse problem (see, e.g., \cite{Gao et al 2014},
\cite{Guliyev}, \cite{Kammanee Bockman 2009}, \cite{Savchuk}) is also a
special case of Problem A. Indeed, let $\mu_{j}^{2}$ be the eigenvalues of the
Sturm--Liouville problem for (\ref{SL equation}) with the boundary conditions
\begin{equation}
y^{\prime}(0)-h_{1}y(0)=0\quad\text{and}\quad y(L)=0, \label{boundary1}%
\end{equation}
while $\nu_{j}^{2}$ be the eigenvalues of the Sturm--Liouville problem for
(\ref{SL equation}) with the conditions
\begin{equation}
y^{\prime}(0)-h_{2}y(0)=0\quad\text{and}\quad y(L)=0, \label{boundary2}%
\end{equation}
$h_{1,2}\in\mathbb{C}$, $h_{1}\neq h_{2}$. Consider a solution $u(\rho,x)$ of
(\ref{SL equation}) such that
\[
u(\mu_{j},0)=u(\nu_{j},0)=1\quad\text{for all}\quad j=0,1,\ldots,
\]
and
\[
u^{\prime}(\mu_{j},0)=h_{1},\ u^{\prime}(\nu_{j},0)=h_{2}\quad\text{for
all}\quad j=0,1,\ldots.
\]
Thus, $u(\rho,x)$ is an eigenfunction of the problem (\ref{SL equation}),
(\ref{boundary1}) when $\rho=\mu_{j}$, and that of the problem
(\ref{SL equation}), (\ref{boundary2}) when $\rho=\nu_{j}$. The set of
$\rho_{k}$ can be chosen in the form of $\left\{  \mu_{0},\nu_{0},\mu_{1}%
,\nu_{1},\cdots\right\}  $ (here, in fact, the order does not matter). The set
of the data (\ref{inpdata}) is
\begin{equation}
a_{k}=1\quad\text{for all }k=1,2,\cdots,\ b_{k}=\left\{
\begin{tabular}
[c]{ll}%
$h_{1},$ & when $\rho_{k}=\mu_{j},$\\
$h_{2},$ & when $\rho_{k}=\nu_{j},$%
\end{tabular}
\ \ \right.  \quad\text{and}\quad\ell_{k}=0, \label{Data two-spectra}%
\end{equation}
$k=1,2,\ldots$. The last equality in (\ref{Data two-spectra}) is due to the
fact that the solution $u(\rho_{k},x)$ being an eigenfunction of
(\ref{SL equation}), (\ref{boundary1}) or (\ref{SL equation}),
(\ref{boundary2}) satisfies the homogeneous Dirichlet condition at $x=L$.
Thus, the inverse two-spectra problem is a special case of Problem A.
Analogously, if $\mu_{j}^{2}$ are the eigenvalues of (\ref{SL equation}),
(\ref{boundary1}) and $\nu_{j}^{2}$ are the Dirichlet-Dirichlet eigenvalues of
(\ref{SL equation}) (i.e., $y(0)=y(L)=0$), then for the set of $\rho_{k}$
chosen as before we have
\begin{equation}
a_{k}=1,\quad b_{k}=h_{1}\text{, when }\rho_{k}=\mu_{j}, \label{dataPr1}%
\end{equation}
\begin{equation}
a_{k}=0,\quad b_{k}=1\text{, when }\rho_{k}=\nu_{j}, \label{dataPr2}%
\end{equation}
and $\ell_{k}=0$. Hence this inverse two-spectra problem is also a special
case of Problem A.

Consider the scattering problem for a potential $q(x)$ with a finite support
in $\left[  0,L\right]  $. A plane wave $e^{-i\rho x}$ incoming from $-\infty$
interacts with $q(x)$. A part of it is reflected, while the other part is
transmitted, so that outside the segment $\left[  0,L\right]  $ the resulting
wave has the form
\[
u(\rho,x)=e^{-i\rho x}+R(\rho)e^{i\rho x},\quad x<0
\]
and
\[
u(\rho,x)=T(\rho)e^{-i\rho x},\quad x>L,
\]
where $R(\rho)$ and $T(\rho)$ denote the reflection and transmission
coefficients, respectively. Now, assume that on a set of points $\rho_{k}%
\in\mathbb{C}$ the values of $R(\rho_{k})$ and $T(\rho_{k})$ are known, and
$q(x)$ needs to be recovered. For the corresponding solution $u(\rho_{k},x)$
of (\ref{SL equation}) we have the initial conditions
\[
u(\rho_{k},0)=1+R(\rho_{k})=:a_{k},
\]%
\[
u^{\prime}(\rho_{k},0)=-i\rho_{k}(1-R(\rho_{k}))=:b_{k}%
\]
and additionally,%
\[
\ell_{k}=T(\rho_{k})e^{-i\rho_{k}L}.
\]
Thus, the inverse scattering problem of the recovery of $q(x)$ from the
reflection and transmission coefficients is also a special case of Problem A.

Our approach to solving Problem A is based on two new developments. Most
importantly, in the present work we develop a new procedure for solving
Problem A based on simple relations between solutions of (\ref{SL equation})
satisfying certain prescribed initial conditions at the opposite endpoints. It
allows us to convert the information on the characteristic functions of
certain regular Sturm-Liouville problems into the information on $q(x)$ on the
whole interval. This procedure works equally well both for real and for
complex-valued potentials. The numerical solution of Problem A reduces to a
successive solution of two systems of linear algebraic equations. Additionally
and to the difference from the previous work \cite{KrSpectrumCompl} and
\cite{KKC2022Mathematics}, we use new Neumann series of Bessel functions
(NSBF) representations for solutions of (\ref{SL equation}), which are
obtained with the aid of a result from \cite{KT2017JMP}. Certain disadvantage
of the previously used NSBF representations (from \cite{KNT}) consists in a
necessity of a double numerical differentiation, which is required in a last
step for recovering $q(x)$. The new NSBF representations allow us to recover
$q(x)$ without any differentiation. The potential is recovered from an
arithmetic combination of the first two coefficients of the representation.
This considerably improves the final accuracy of the result. However, the use
of the new NSBF representations is restricted to smooth potentials, $q\in
C^{1}\left[  0,L\right]  $, that reduces their applicability. Thus, in the
present work, a new method for solving Problem A for complex-valued
continuously differentiable potentials is developed. It is easy for
implementation, direct, accurate and applicable to a large variety of inverse problems.

Besides this introduction the paper contains six sections. In Section
\ref{Sect Prelim} we present the NSBF representations for solutions of
(\ref{SL equation}). Here Theorem \ref{Th NSBF} is not new and can be found in
\cite{KNT}, and Theorem \ref{Th NSBF new} together with Remark
\ref{Rem final form of NSBF} present the new NSBF representations, which are
obtained from those of \cite{KT2017JMP}. In Section \ref{Sect comp phi and S}
we explain, how the characteristic functions of two regular Sturm-Liouville
problems can be obtained directly from the input data of the form
(\ref{inpdata}), as well as some other relevant parameters. In Section
\ref{Sect Recovery q} we show, how these data lead to the recovery of $q(x)$.
In Section \ref{Sect Summary} we summarize the method. Section
\ref{Sect Numeric} contains numerical results and discussion. Finally, Section
\ref{Sect Conclusions} presents some concluding remarks.

\section{Neumann series of Bessel functions representations for solutions
\label{Sect Prelim}}

By $\varphi(\rho,x)$ and $S(\rho,x)$ we denote the solutions of the equation
\begin{equation}
-y^{\prime\prime}(x)+q(x)y(x)=\rho^{2}y(x),\quad x\in(0,L) \label{Schr}%
\end{equation}
satisfying the initial conditions
\[
\varphi(\rho,0)=1,\quad\varphi^{\prime}(\rho,0)=0,
\]%
\[
S(\rho,0)=0,\quad S^{\prime}(\rho,0)=1.
\]
In \cite{KNT} (see also \cite[Sect. 9.2]{KrBook2020}) convenient series
representations for the solutions $\varphi(\rho,x)$ and $S(\rho,x)$ of
(\ref{Schr}) were obtained.

\begin{theorem}
[\cite{KNT}]\label{Th NSBF} Let $q\in%
\mathcal{L}%
_{2}(0,L)$. The solutions $\varphi(\rho,x)$ and $S(\rho,x)$ admit the
following series representations
\begin{align}
\varphi(\rho,x)  &  =\cos\left(  \rho x\right)  +\sum_{n=0}^{\infty}%
(-1)^{n}g_{n}(x)\mathbf{j}_{2n}(\rho x),\label{phiNSBF}\\
S(\rho,x)  &  =\frac{\sin\left(  \rho x\right)  }{\rho}+\frac{1}{\rho}%
\sum_{n=0}^{\infty}(-1)^{n}s_{n}(x)\mathbf{j}_{2n+1}(\rho x), \label{S}%
\end{align}
where $\mathbf{j}_{k}(z)$ stands for the spherical Bessel function of order
$k$ (see, e.g., \cite{AbramowitzStegunSpF}). The coefficients $g_{n}(x)$ and
$s_{n}(x)$ can be calculated following a simple recurrent integration
procedure (see \cite{KNT} or \cite[Sect. 9.4]{KrBook2020}), starting with
\begin{equation}
g_{0}(x)=\varphi(0,x)-1\quad\text{and}\quad s_{0}(x)=3\left(  \frac{S(0,x)}%
{x}-1\right)  . \label{beta0}%
\end{equation}
For every $\rho\in\mathbb{C}$ the series converge pointwise. For every
$x\in\left[  0,L\right]  $ the series converge uniformly on any compact set of
the complex plane of the variable $\rho$, and the remainders of their partial
sums admit estimates independent of $\operatorname{Re}\rho$.
\end{theorem}

This last feature of the series representations (the independence of the
bounds for the remainders of $\operatorname{Re}\rho$) is a direct consequence
of the fact that the representations are obtained by expanding the integral
kernels of the transmutation operators (for their theory we refer to
\cite{LevitanInverse}, \cite{Marchenko}, \cite{SitnikShishkina Elsevier}) into
Fourier-Legendre series (see \cite{KNT} and \cite[Sect. 9.2]{KrBook2020}). In
particular, it means that for $\varphi_{N}(\rho,x):=\cos\left(  \rho x\right)
+\sum_{n=0}^{N}(-1)^{n}g_{n}(x)\mathbf{j}_{2n}(\rho x)$ and $S_{N}%
(\rho,x):=\frac{\sin\left(  \rho x\right)  }{\rho}+\frac{1}{\rho}\sum
_{n=0}^{N}(-1)^{n}s_{n}(x)\mathbf{j}_{2n+1}(\rho x)$ the estimates hold
\begin{equation}
\left\vert \varphi(\rho,x)-\varphi_{N}(\rho,x)\right\vert \leq\frac
{2\varepsilon_{N}(x)\,\sinh(Cx)}{C}\quad\text{and}\quad\left\vert
S(\rho,x)-S_{N}(\rho,x)\right\vert \leq\frac{2\varepsilon_{N}(x)\,\sinh
(Cx)}{C} \label{estc2}%
\end{equation}
for any $\rho\in\mathbb{C}$ belonging to the strip $\left\vert
\operatorname{Im}\rho\right\vert \leq C$, $C\geq0$, where $\varepsilon_{N}(x)$
is a positive function tending to zero when $N\rightarrow\infty$. That is,
roughly speaking, the approximate solutions $\varphi_{N}(\rho,x)$ and
$S_{N}(\rho,x)$ approximate the exact ones equally well for small and for
large values of $\operatorname{Re}\rho$. Detailed estimates for the series
remainders depending on the regularity of the potential can be found in
\cite{KNT}.

Note that formulas (\ref{beta0}) indicate that the potential $q(x)$ can be
recovered from the first coefficients of the series (\ref{phiNSBF}) or
(\ref{S}). We have
\begin{equation}
q(x)=\frac{g_{0}^{\prime\prime}(x)}{g_{0}(x)+1}\quad\text{and\quad}%
q(x)=\frac{\left(  xs_{0}(x)\right)  ^{\prime\prime}}{xs_{0}(x)+3x}.
\label{qi from g0}%
\end{equation}

We recall that any series of the type $\sum_{n=0}^{\infty}a_{n}J_{\nu+n}(z)$
is called a Neumann series of Bessel functions (NSBF) (see \cite{Baricz et
al}, \cite{Baricz et al Book}, \cite[Chapter XVI]{Watson}, \cite{Wilkins}).

Together with the solutions $\varphi(\rho,x)$ and $S(\rho,x)$ we will also
need to consider a solution $T(\rho,x)$ of (\ref{Schr}) which satisfies the
initial conditions at $x=L$:%
\begin{equation}
T(\rho,L)=0\quad\text{and}\quad T^{\prime}(\rho,L)=1. \label{init T}%
\end{equation}

Note that
\begin{equation}
T(\rho,x)=\varphi(\rho,L)S(\rho,x)-S(\rho,L)\varphi(\rho,x). \label{solT}%
\end{equation}
Indeed, the first condition in (\ref{init T}) is obviously fulfilled by the
function (\ref{solT}), and to verify the second one, from (\ref{solT}) we have%
\[
T^{\prime}(\rho,L)=\varphi(\rho,L)S^{\prime}(\rho,L)-S(\rho,L)\varphi^{\prime
}(\rho,L)=W\left[  \varphi,S\right]  \left(  L\right)  ,
\]
where $W\left[  \varphi,S\right]  \left(  L\right)  $ stands for the Wronskian
of $\varphi(\rho,x)$ and $S(\rho,x)$ at the point $x=L$. Since the Wronskian
of two solutions of (\ref{Schr}) is constant (see, e.g., \cite[p.
327]{Hartman}), we have that
\[
W\left[  \varphi,S\right]  \left(  L\right)  =W\left[  \varphi,S\right]
\left(  0\right)  =\varphi(\rho,0)S^{\prime}(\rho,0)-S(\rho,0)\varphi^{\prime
}(\rho,0)=1,
\]
and thus the solution (\ref{solT}) satisfies both initial conditions
(\ref{init T}).

The approach developed for solving inverse problems is based on the NSBF
representations and, in particular, on the possibility of recovering the
potential from the first coefficients. However, the use of the representations
from Theorem \ref{Th NSBF} requires a double differentiation in the final step
(formulas (\ref{qi from g0})). Although, as it was shown in a number of
publications \cite{Kr2019JIIP}, \cite{KrBook2020}, \cite{KrSpectrumCompl},
\cite{KKC2022Mathematics}, \cite{KST2020IP}, \cite{KT2021 IP1}, \cite{KT2021
IP2}, this is not a serious drawback and leads to a satisfactory final
accuracy of the recovered potential, in the present work another possibility
is explored, which does not require any differentiation in this last step and
results in a still more accurate recovery of the potential. Namely, we use the
following NSBF representations for the solutions $\varphi(\rho,x)$,
$S(\rho,x)$ and $T(\rho,x)$.

\begin{theorem}
\label{Th NSBF new} Let $q\in C^{1}\left[  0,L\right]  $. Then the solutions
$\varphi(\rho,x)$, $S(\rho,x)$ and $T(\rho,x)$ of (\ref{Schr}) admit the
series representations%
\begin{equation}
\varphi(\rho,x)=\cos\left(  \rho x\right)  +\frac{\sin\left(  \rho x\right)
}{\rho}\omega(x)+\frac{\cos\left(  \rho x\right)  }{\rho^{2}}q^{-}(x)-\frac
{1}{\rho^{2}}\sum_{n=0}^{\infty}(-1)^{n}\varphi_{n}(x)\mathbf{j}_{2n}(\rho x),
\label{phi NSBF new}%
\end{equation}%
\begin{equation}
S(\rho,x)=\frac{\sin\left(  \rho x\right)  }{\rho}-\frac{\cos\left(  \rho
x\right)  }{\rho^{2}}\omega(x)+\frac{\sin\left(  \rho x\right)  }{\rho^{3}%
}q^{+}(x)-\frac{1}{\rho^{3}}\sum_{n=0}^{\infty}(-1)^{n}\sigma_{n}%
(x)\mathbf{j}_{2n+1}(\rho x), \label{S NSBF new}%
\end{equation}%
\begin{align}
T(\rho,x)  &  =-\frac{\sin\left(  \rho\left(  L-x\right)  \right)  }{\rho
}+\frac{\cos\left(  \rho\left(  L-x\right)  \right)  }{\rho^{2}}\omega
_{L}(x)-\frac{\sin\left(  \rho\left(  L-x\right)  \right)  }{\rho^{3}}%
q_{L}^{+}(x)\nonumber\\
&  +\frac{1}{\rho^{3}}\sum_{n=0}^{\infty}(-1)^{n}\theta_{n}(x)\mathbf{j}%
_{2n+1}(\rho\left(  L-x\right)  ), \label{T NSBF new}%
\end{align}
where
\[
\omega(x):=\frac{1}{2}\int_{0}^{x}q(s)ds,\quad\omega_{L}(x):=\frac{1}{2}%
\int_{x}^{L}q(s)ds,
\]%
\[
q^{\pm}(x):=\frac{q(x)\pm q(0)}{4}-\frac{\omega^{2}(x)}{2},\quad q_{L}%
^{+}(x):=\frac{q(x)+q(L)}{4}-\frac{\omega_{L}^{2}(x)}{2}.
\]
For every $x\in\left[  0,L\right]  $ the series converge uniformly on any
compact set of the complex plane of the variable $\rho$, and for any $\rho
\in\mathbb{C}\backslash\left\{  0\right\}  $ the remainders of their partial
sums admit the estimates%
\[
\left\vert \varphi(\rho,x)-\varphi_{N}(\rho,x)\right\vert \leq\frac
{\varepsilon_{N}(x)\,}{\left\vert \rho\right\vert ^{2}}\sqrt{\frac
{\sinh(2\operatorname{Im}\rho x)}{\operatorname{Im}\rho}},
\]%
\[
\left\vert S(\rho,x)-S_{N}(\rho,x)\right\vert \leq\frac{\varepsilon_{N}%
(x)\,}{\left\vert \rho\right\vert ^{3}}\sqrt{\frac{\sinh(2\operatorname{Im}%
\rho x)}{\operatorname{Im}\rho}},
\]%
\begin{equation}
\left\vert T(\rho,x)-T_{N}(\rho,x)\right\vert \leq\frac{\varepsilon_{N}%
(x)\,}{\left\vert \rho\right\vert ^{3}}\sqrt{\frac{\sinh(2\operatorname{Im}%
\rho x)}{\operatorname{Im}\rho}}, \label{estimTnew}%
\end{equation}
where the subindex $N$ indicates that in (\ref{phi NSBF new}%
)-(\ref{T NSBF new}) the sum is taken up to $N$, and $\varepsilon_{N}(x)$ is a
positive function tending to zero when $N\rightarrow\infty$.
\end{theorem}

\emph{Proof.}\textbf{ }This theorem is a direct corollary of Theorem 10 from
\cite{KT2017JMP}. The representations (\ref{phi NSBF new}) and
(\ref{S NSBF new}) are obtained from \cite[formula (23)]{KT2017JMP} by
considering the relations $\varphi(\rho,x)=\left(  u(\rho,x)+u(-\rho
,x)\right)  /2$ and $S(\rho,x)=\left(  u(\rho,x)-u(-\rho,x)\right)  /(2i\rho)$
with $u(\rho,x)$ being the solution considered in \cite[formula (23)]%
{KT2017JMP}, which satisfies the conditions $u(\rho,0)=1$ and $u^{\prime}%
(\rho,0)=i\rho$. The representation (\ref{T NSBF new}) is obtained from
(\ref{S NSBF new}) by flipping the interval. The estimates follow from
\cite[formula (24)]{KT2017JMP}. $\blacksquare$

\begin{remark}
\label{Rem final form of NSBF}In \cite{KT2017JMP} it was shown that
\[
\varphi_{0}(x)=q^{-}(x)
\]
and
\begin{equation}
\sigma_{0}(x)=3\left(  q^{+}(x)-\frac{\omega(x)}{x}\right)  . \label{sigma0}%
\end{equation}
Thus, the series representations (\ref{phi NSBF new}) and (\ref{S NSBF new})
can be written in a more convenient form by regrouping the terms. We have that
the coefficient at $\varphi_{0}(x)$ in (\ref{phi NSBF new}) can be written as%
\[
\frac{\cos\left(  \rho x\right)  }{\rho^{2}}-\frac{\mathbf{j}_{0}(\rho
x)}{\rho^{2}}=\frac{1}{\rho^{2}}\left(  \cos\left(  \rho x\right)  -\frac
{\sin\left(  \rho x\right)  }{\rho x}\right)  =-\frac{x\mathbf{j}_{1}(\rho
x)}{\rho},
\]
where we used the explicit forms of the spherical Bessel functions
\[
\mathbf{j}_{0}(z)=\frac{\sin z}{z}\quad\text{and}\quad\mathbf{j}_{1}%
(z)=\frac{\sin z}{z^{2}}-\frac{\cos z}{z}.
\]
Thus,
\begin{equation}
\varphi(\rho,x)=\cos\left(  \rho x\right)  +\frac{\sin\left(  \rho x\right)
}{\rho}\omega(x)-\frac{x\mathbf{j}_{1}(\rho x)}{\rho}q^{-}(x)-\frac{1}%
{\rho^{2}}\sum_{n=1}^{\infty}(-1)^{n}\varphi_{n}(x)\mathbf{j}_{2n}(\rho x).
\label{phi1}%
\end{equation}

Similarly, taking into account (\ref{sigma0}), we can regroup the second,
third and fourth terms in (\ref{S NSBF new}):%
\begin{align*}
&  -\frac{\cos\left(  \rho x\right)  }{\rho^{2}}\omega(x)+\frac{\sin\left(
\rho x\right)  }{\rho^{3}}q^{+}(x)-\frac{3\mathbf{j}_{1}(\rho x)}{\rho^{3}%
}\left(  q^{+}(x)-\frac{\omega(x)}{x}\right) \\
&  =\omega(x)\left(  \frac{3\mathbf{j}_{1}(\rho x)}{\rho^{3}x}-\frac
{\cos\left(  \rho x\right)  }{\rho^{2}}\right)  +q^{+}(x)\left(  \frac
{\sin\left(  \rho x\right)  }{\rho^{3}}-\frac{3\mathbf{j}_{1}(\rho x)}%
{\rho^{3}}\right)  .
\end{align*}
Thus,
\begin{align}
S(\rho,x)  &  =\frac{\sin\left(  \rho x\right)  }{\rho}+\frac{\omega(x)}%
{\rho^{2}}\left(  \frac{3\mathbf{j}_{1}(\rho x)}{\rho x}-\cos\left(  \rho
x\right)  \right)  +\frac{q^{+}(x)}{\rho^{3}}\left(  \sin\left(  \rho
x\right)  -3\mathbf{j}_{1}(\rho x)\right) \nonumber\\
&  -\frac{1}{\rho^{3}}\sum_{n=1}^{\infty}(-1)^{n}\sigma_{n}(x)\mathbf{j}%
_{2n+1}(\rho x). \label{S1}%
\end{align}

Analogously,
\begin{align}
T(\rho,x)  &  =-\frac{\sin\left(  \rho\left(  L-x\right)  \right)  }{\rho
}-\frac{\omega_{L}(x)}{\rho^{2}}\left(  \frac{3\mathbf{j}_{1}(\rho\left(
L-x\right)  )}{\rho\left(  L-x\right)  }-\cos\left(  \rho\left(  L-x\right)
\right)  \right) \nonumber\\
&  -\frac{q_{L}^{+}(x)}{\rho^{3}}\left(  \sin\left(  \rho\left(  L-x\right)
\right)  -3\mathbf{j}_{1}(\rho\left(  L-x\right)  )\right)  +\frac{1}{\rho
^{3}}\sum_{n=1}^{\infty}(-1)^{n}\theta_{n}(x)\mathbf{j}_{2n+1}(\rho\left(
L-x\right)  ). \label{T1}%
\end{align}

\end{remark}

\begin{remark}
Unlike the representations from Theorem \ref{Th NSBF}, which are applicable in
a general case $q\in%
\mathcal{L}%
_{2}(0,L)$ (and in fact even for $q$ from a larger class $W_{2}^{-1}(0,L)$,
see \cite{KV2023arxiv} and references therein), the representations
(\ref{phi1}), (\ref{S1}) and (\ref{T1}) are applicable to $q\in C^{1}\left[
0,L\right]  $. However, in this case the potential can be recovered without
any differentiation in the last step, as we explain further, in subsection
\ref{subsect recovery q}, that leads to a higher accuracy in recovering $q(x)$.
\end{remark}

Thus, in the present work, for solving inverse coefficient problems we use the
NSBF representations (\ref{phi1})-(\ref{T1}), while for computing the input
data (direct problem) we use the NSBF representations from Theorem
\ref{Th NSBF}.

\section{Computation of $\varphi(\rho,L)$ and $S(\rho,L)$ from the input data
(\ref{inpdata})\label{Sect comp phi and S}}

Here we adapt the procedure from \cite{KKC2022Mathematics} to the new NSBF
representations (\ref{phi1}), (\ref{S1}). Denote $a_{k}:=a(\rho_{k})$ and
$b_{k}:=b(\rho_{k})$. We have $u(\rho_{k},x)=a_{k}\varphi(\rho_{k}%
,x)+b_{k}S(\rho_{k},x)$, and hence
\[
a_{k}\varphi(\rho_{k},L)+b_{k}S(\rho_{k},L)=\ell_{k}.
\]
Substitution of (\ref{phi1}) and (\ref{S1}) into this equality gives us the
equation%
\begin{gather*}
\left(  \frac{a_{k}\sin\left(  \rho_{k}L\right)  }{\rho_{k}}+\frac{b_{k}}%
{\rho_{k}^{2}}\left(  \frac{3\mathbf{j}_{1}(\rho_{k}L)}{\rho_{k}L}-\cos\left(
\rho_{k}L\right)  \right)  \right)  \omega(L)-\frac{a_{k}L\mathbf{j}_{1}%
(\rho_{k}L)}{\rho_{k}}q^{-}(L)-\frac{a_{k}}{\rho_{k}^{2}}\sum_{n=1}^{\infty
}(-1)^{n}\varphi_{n}(L)\mathbf{j}_{2n}(\rho_{k}L)\\
+\frac{b_{k}}{\rho_{k}^{3}}\left(  \sin\left(  \rho_{k}L\right)
-3\mathbf{j}_{1}(\rho_{k}L)\right)  q^{+}(L)-\frac{b_{k}}{\rho_{k}^{3}}%
\sum_{n=1}^{\infty}(-1)^{n}\sigma_{n}(L)\mathbf{j}_{2n+1}(\rho_{k}L)=\ell
_{k}-a_{k}\cos\left(  \rho_{k}L\right)  -\frac{b_{k}}{\rho_{k}}\sin\left(
\rho_{k}L\right)
\end{gather*}
for all $\rho_{k}$. This leads to a finite system of linear algebraic
equations for computing the NSBF coefficients $\omega(L)$, $q^{-}(L)$,
$q^{+}(L)$, $\varphi_{n}(L)$, $\sigma_{n}(L)$, $n=1,\ldots,N$,
\begin{gather}
\left(  \frac{a_{k}\sin\left(  \rho_{k}L\right)  }{\rho_{k}}+\frac{b_{k}}%
{\rho_{k}^{2}}\left(  \frac{3\mathbf{j}_{1}(\rho_{k}L)}{\rho_{k}L}-\cos\left(
\rho_{k}L\right)  \right)  \right)  \omega(L)-\frac{a_{k}L\mathbf{j}_{1}%
(\rho_{k}L)}{\rho_{k}}q^{-}(L)-\frac{a_{k}}{\rho_{k}^{2}}\sum_{n=1}%
^{N}(-1)^{n}\varphi_{n}(L)\mathbf{j}_{2n}(\rho_{k}L)\nonumber\\
+\frac{b_{k}}{\rho_{k}^{3}}\left(  \sin\left(  \rho_{k}L\right)
-3\mathbf{j}_{1}(\rho_{k}L)\right)  q^{+}(L)-\frac{b_{k}}{\rho_{k}^{3}}%
\sum_{n=1}^{N}(-1)^{n}\sigma_{n}(L)\mathbf{j}_{2n+1}(\rho_{k}L)=\ell_{k}%
-a_{k}\cos\left(  \rho_{k}L\right)  -\frac{b_{k}}{\rho_{k}}\sin\left(
\rho_{k}L\right)  \label{FirstSystem}%
\end{gather}
for $k=1,\ldots,K$.

\begin{remark}
Here and below we do not look to deal with square systems of equations. In
computations a least-squares solution of an overdetermined system gives
satisfactory results and allows us to make use of all available data, while
keeping the number of the coefficients relatively small (in practice $N=4\ $or
$5$ may result sufficient). Thus, $K\geq2N+3$.
\end{remark}

\begin{remark}
\label{Rem omega and all that}The value $\omega(L)$ is an important parameter,
that arises in the second term of the asymptotics of the eigenvalues of the
Sturm-Liouville problems for (\ref{Schr}). Many available in bibliography
numerical techniques for solving inverse Sturm-Liouville problems require its
prior knowledge as, e.g., \cite{Rundell Sacks} or \cite{IgnatievYurko}. Here
we obtain it immediately, directly from the input data of the problem, as well
as the parameters $q^{-}(L)$ and $q^{+}(L)$. This gives us additionally the
values of the potential $q(x)$ at the end points:%
\[
q(0)=2\left(  q^{+}(L)-q^{-}(L)\right)  \quad\text{and\quad}q(L)=2\left(
q^{+}(L)+q^{-}(L)+\omega^{2}(L)\right)  .
\]

\end{remark}

The knowledge of the coefficients $\omega(L)$, $q^{-}(L)$, $q^{+}(L)$,
$\varphi_{n}(L)$, $\sigma_{n}(L)$, $n=1,\ldots,N$ allows us to compute the
functions $\varphi_{N}(\rho,L)$ and $S_{N}(\rho,L)$ for any value of $\rho$.
Due to the estimates from Theorem \ref{Th NSBF new} we know that the accuracy
of the approximation of the exact solutions $\varphi(\rho,L)$ and $S(\rho,L)$
by these approximate ones does not deteriorate for large values of $\rho
\in\mathbb{R}$ and even improves. So, in fact, we deal with the problem of
converting the knowledge of $\varphi(\rho,L)$ and $S(\rho,L)$ into the
knowledge of $q(x)$ on the whole interval.

\section{Recovery of $q(x)$ from $\varphi(\rho,L)$ and $S(\rho,L)$
\label{Sect Recovery q}}

\subsection{Main system of linear algebraic equations}

Consider the identity (\ref{solT}), in which $\varphi(\rho,L)$ and $S(\rho,L)$
are already known. We use (\ref{solT}) to construct the main system of linear
algebraic equations. Assume $\varphi(\rho,L)$ and $S(\rho,L)$ to be computed
on a countable set of points $\left\{  \gamma_{k}^{2}\right\}  _{k=1}^{\infty
}$. Denote
\begin{equation}
S_{k}:=S(\gamma_{k},L)\quad\text{and\quad}F_{k}:=\varphi(\gamma_{k}%
,L).\label{T_k}%
\end{equation}
Substitution of the NSBF representations (\ref{phi1}), (\ref{S1}) and
(\ref{T1}) into (\ref{solT}) leads to a system of linear algebraic equations
for the NSBF coefficients.

\begin{theorem}
Let $q\in C^{1}\left[  0,L\right]  $. Then for all $x\in\left(  0,L\right)  $
the functions
\[
\omega(x),\qquad Q(x):=\frac{q(x)}{4}-\frac{\omega^{2}(x)}{2},\qquad
q_{0}:=\frac{q(0)}{4},
\]
$\omega_{L}(x)$, $q_{L}^{+}(x)$ and $\left\{  \varphi_{n}(x),\,\sigma
_{n}(x),\,\theta_{n}(x)\right\}  _{n=1}^{\infty}$ satisfy the system of linear
algebraic equations%
\[
A_{k1}(x)\omega(x)+A_{k2}(x)Q(x)+A_{k3}(x)q_{0}+A_{k4}(x)\omega_{L}%
(x)+A_{k5}(x)q_{L}^{+}(x)
\]%
\[
+\sum_{n=1}^{\infty}B_{kn}(x)\varphi_{n}(x)+\sum_{n=1}^{\infty}C_{kn}%
(x)\sigma_{n}(x)+\sum_{n=1}^{\infty}D_{kn}(x)\theta_{n}(x)
\]%
\begin{equation}
=-\frac{\sin\left(  \gamma_{k}\left(  L-x\right)  \right)  }{\gamma_{k}}%
+S_{k}\cos\left(  \gamma_{k}x\right)  -\frac{F_{k}\sin\left(  \gamma
_{k}x\right)  }{\gamma_{k}},\quad k=1,2,\ldots,\label{main system}%
\end{equation}
where%
\[
A_{k1}(x):=-\frac{S_{k}\sin\left(  \gamma_{k}x\right)  }{\gamma_{k}}%
+\frac{F_{k}}{\gamma_{k}^{2}}\left(  \frac{3\mathbf{j}_{1}(\gamma_{k}%
x)}{\gamma_{k}x}-\cos\left(  \gamma_{k}x\right)  \right)  ,
\]%
\[
A_{k2}(x):=\frac{F_{k}}{\gamma_{k}^{3}}\left(  \sin\left(  \gamma_{k}x\right)
-3\mathbf{j}_{1}(\gamma_{k}x)\right)  +\frac{S_{k}x\mathbf{j}_{1}\left(
\gamma_{k}x\right)  }{\gamma_{k}},
\]%
\[
A_{k3}(x):=-\frac{S_{k}x\mathbf{j}_{1}\left(  \gamma_{k}x\right)  }{\gamma
_{k}}+\frac{F_{k}}{\gamma_{k}^{3}}\left(  \sin\left(  \gamma_{k}x\right)
-3\mathbf{j}_{1}(\gamma_{k}x)\right)  ,
\]%
\[
A_{k4}(x):=\frac{1}{\gamma_{k}^{2}}\left(  \frac{3\mathbf{j}_{1}(\gamma
_{k}\left(  L-x\right)  )}{\gamma_{k}\left(  L-x\right)  }-\cos\left(
\gamma_{k}\left(  L-x\right)  \right)  \right)  ,
\]%
\[
A_{k5}(x):=\frac{1}{\gamma_{k}^{3}}\left(  \sin\left(  \gamma_{k}\left(
L-x\right)  \right)  -3\mathbf{j}_{1}(\gamma_{k}\left(  L-x\right)  )\right)
,
\]%
\begin{align*}
B_{kn}(x) &  :=\left(  -1\right)  ^{n}\frac{S_{k}}{\gamma_{k}^{2}}%
\mathbf{j}_{2n}(\gamma_{k}x),\quad C_{kn}(x):=\left(  -1\right)  ^{n+1}%
\frac{F_{k}}{\gamma_{k}^{3}}\mathbf{j}_{2n+1}(\gamma_{k}x),\\
D_{kn}(x) &  :=\left(  -1\right)  ^{n+1}\frac{1}{\gamma_{k}^{3}}%
\mathbf{j}_{2n+1}(\gamma_{k}\left(  L-x\right)  ).
\end{align*}

\end{theorem}

\emph{Proof.}\textbf{ }The proof of (\ref{main system}) consists in
substituting the series representations (\ref{phi1})-(\ref{T1}) into
(\ref{solT}) and grouping coefficients at $\omega(x)$, $Q(x)$ and $q_{0}$.
$\blacksquare$

\begin{remark}
Note that $q_{0}=q(0)/4$ is computed in the first step, when solving system
(\ref{FirstSystem}), see Remark \ref{Rem omega and all that}. Thus, in the
system (\ref{main system}) the term $A_{k3}(x)q_{0}$ can be treated as known
and moved to the right-hand side. Moreover, we have the following relations%
\[
\omega_{L}(x)=\omega(L)-\omega(x)
\]
and
\[
q_{L}^{+}(x)=\frac{q(x)+q(L)}{4}-\frac{1}{2}\left(  \omega(L)-\omega
(x)\right)  ^{2}=\omega(L)\omega(x)+Q(x)+\frac{q(L)}{4}-\frac{\omega^{2}%
(L)}{2}.
\]
The values of $\omega(L)$ and $q(L)$ are also computed in the first step, see
Remark \ref{Rem omega and all that}. Thus, the corresponding terms in
(\ref{main system}) can also be moved to the right-hand side, and
(\ref{main system}) takes the form%
\begin{align*}
& \left(  A_{k1}(x)-A_{k4}(x)+\omega(L)A_{k5}(x)\right)  \omega(x)+\left(
A_{k2}(x)+A_{k5}(x)\right)  Q(x)\\
& +\sum_{n=1}^{\infty}B_{kn}(x)\varphi_{n}(x)+\sum_{n=1}^{\infty}%
C_{kn}(x)\sigma_{n}(x)+\sum_{n=1}^{\infty}D_{kn}(x)\theta_{n}(x)
\end{align*}%
\begin{gather*}
=-\frac{\sin\left(  \gamma_{k}\left(  L-x\right)  \right)  }{\gamma_{k}}%
+S_{k}\cos\left(  \gamma_{k}x\right)  -\frac{F_{k}\sin\left(  \gamma
_{k}x\right)  }{\gamma_{k}}\\
-A_{k3}(x)q_{0}-A_{k4}(x)\omega(L)+\frac{A_{k5}(x)}{2}\left(  \omega
^{2}(L)-\frac{q(L)}{2}\right)  ,\quad k=1,2,\ldots.
\end{gather*}

\end{remark}

\subsection{Final step: recovery of $q(x)$\label{subsect recovery q}}

From the data obtained after solving system (\ref{main system}) the potential
$q(x)$ can be recovered in the following ways.

1. Differentiation of $2\omega(x)$ gives us $q(x)$. Moreover, since typically
the numerical differentiation provides less accurate results at the endpoints
of the interval, it is convenient to use the values $q(0)$ and $q(L)$ which
are obtained in the first step by solving\ system (\ref{FirstSystem}), see
Remark \ref{Rem omega and all that}. Thus, the potential $q(x)$ is obtained
from $q(x)=2\omega^{\prime}(x)$, and at the endpoints its previously computed
values are used.

2. The potential $q(x)$ can be obtained without differentiating, directly from
$Q(x)$ and $\omega(x)$, which are computed from (\ref{main system}). We have%
\[
q(x)=4Q(x)+2\omega^{2}(x).
\]
Again, the values $q(0)$ and $q(L)$ are obtained in the first step by
solving\ system (\ref{FirstSystem}), see Remark \ref{Rem omega and all that}.

Both options give similar numerical results.

\section{Summary of the method\label{Sect Summary}}

As it is seen from the previous pages, solution of the considered inverse
coefficient problem reduces to the solution of the main system of linear
algebraic equations (\ref{main system}) in its reduced form
\[
\left(  A_{k1}(x)-A_{k4}(x)+\omega(L)A_{k5}(x)\right)  \omega(x)+\left(
A_{k2}(x)+A_{k5}(x)\right)  Q(x)
\]%
\[
+\sum_{n=1}^{N}B_{kn}(x)\varphi_{n}(x)+\sum_{n=1}^{N}C_{kn}(x)\sigma
_{n}(x)+\sum_{n=1}^{N}D_{kn}(x)\theta_{n}(x)
\]%
\[
=-\frac{\sin\left(  \gamma_{k}\left(  L-x\right)  \right)  }{\gamma_{k}}%
+S_{k}\cos\left(  \gamma_{k}x\right)  -\frac{F_{k}\sin\left(  \gamma
_{k}x\right)  }{\gamma_{k}}-A_{k3}(x)q_{0}%
\]%
\begin{equation}
-A_{k4}(x)\omega(L)+\frac{A_{k5}(x)}{2}\left(  \omega^{2}(L)-\frac{q(L)}%
{2}\right)  ,\quad k=1,\ldots,m.\label{main system reduced}%
\end{equation}
Having solved it, one can compute $q(x)$ easily, as explained in subsection
\ref{subsect recovery q}.

The proposed algorithm for solving Problem A can be summarized as follows.

1. For a given set of input data (\ref{inpdata}) solve (\ref{FirstSystem}).
This gives us $\omega(L)$, $q^{-}(L)$, $q^{+}(L)$, $\left\{  \varphi
_{n}(L),\,\sigma_{n}(L)\right\}  _{n=1}^{N}$. As explained in Remark
\ref{Rem omega and all that}, we obtain $q(0)$ and $q(L)$.

2. Choose a set of points $\left\{  \gamma_{k}\right\}  _{k=1}^{m}$ and
compute (\ref{T_k}).

3. Solve (\ref{main system reduced}).

4. Recover $q(x)$ as explained in subsection \ref{subsect recovery q}.

Note that the choice of the distribution of points $\gamma_{k}$ as well as
their total number are arbitrary. The estimates for the remainders of the NSBF
series involved suggest that it is convenient to choose $\gamma_{k}$ in a
strip $\left\vert \operatorname{Im}\rho\right\vert \leq C$, with $C\geq0$
being not too large. Numerical experiments show that a uniform distribution of
$\gamma_{k}$ can be less advantageous than, e.g., a logarithmically equally
spaced points distribution, with more densely distributed points near the
origin and less densely far from it.

\section{Numerical results\label{Sect Numeric}}

The proposed approach can be implemented directly using an available numerical
computing environment. All the reported computations were performed in Matlab
2017 on an Intel i7-7600U equipped laptop computer and took no more than
several seconds.

\textbf{Example 1 }Let us start with the solution of the inverse two-spectra
problem with the data (\ref{dataPr1}) with $h_{1}=0$ and (\ref{dataPr2}) for
the potential
\[
q(x)=e^{x}+i,\quad x\in\left[  0,\pi\right]  .
\]

The Dirichlet-Dirichlet and Neumann-Dirichlet eigenvalues are obtained from
those for the potential $e^{x}$ by adding the imaginary unit. For computing
the \textquotedblleft exact\textquotedblright\ eigenvalues of $e^{x}$ we used
the Matslise package \cite{Ledoux et al}. In Fig. 1 the real and imaginary
parts of the potential recovered from $10$ eigenpairs are presented. The
maximum absolute error was $0.056$, attained with both ways for the recovery
of $q(x)$ in the last step from subsection \ref{subsect recovery q}. The
absolute error of $\omega_{L}(0)$ computed in the first step was $0.00058$.
The points $\left\{  \gamma_{k}\right\}  _{k=1}^{m}$ in step 2 were chosen as
$\gamma_{k}=10^{\alpha_{k}}$ with $\alpha_{k}$ being distributed uniformly on
$\left[  \lg(0.1),\lg(1500)\right]  \approx\left[  -1,\,3.176\right]  $. This
is a logarithmically equally spaced points distribution on the segment
$\left[  0.1,\,1500\right]  $. The number $m$ was chosen as $m=700$. The
number $N$ in (\ref{main system reduced}) was chosen as $N=12$.%

\begin{figure}
[ptb]
\begin{center}
\includegraphics[
height=3.5572in,
width=4.7386in
]%
{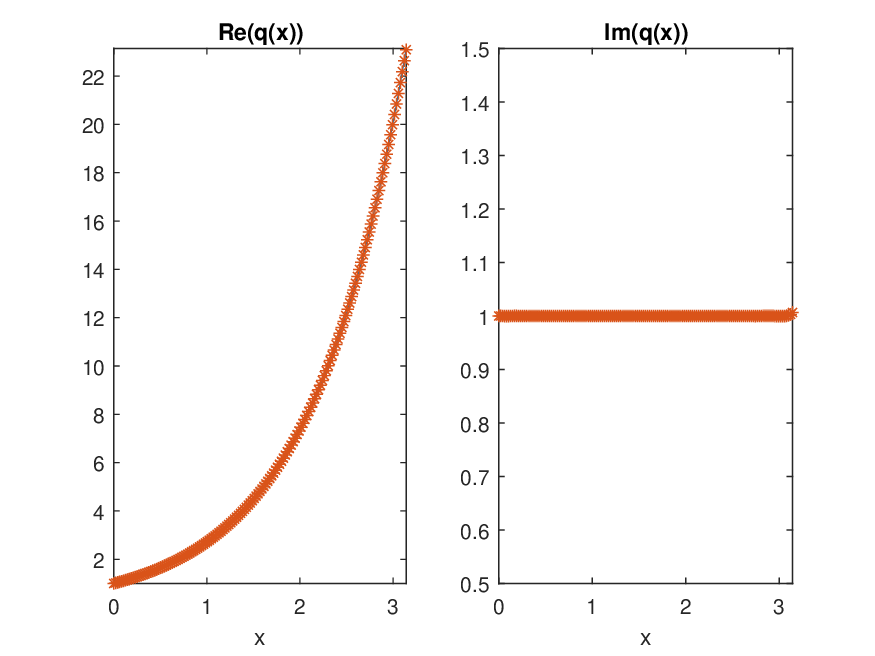}%
\caption{Real and imaginary parts of the potential from Example 1, recovered
from 10 eigenpairs. The maximum absolute error of the recovered potential is
$0.056$.}%
\label{Fig1}%
\end{center}
\end{figure}
The result of the recovery of $q(x)$ from $15$ eigenpairs and preserving the
same values of the rest of the parameters was considerably more precise. The
parameter $\omega_{L}(0)$ computed in the first step was obtained with the
absolute error $2.8\cdot10^{-9}$. The maximum absolute error of $q(x)$
recovered following the first option from subsection \ref{subsect recovery q}
was $1.04\cdot10^{-6}$, while following the second one $0.00024$.

\textbf{Example 2. }Consider Problem A for the potential
\[
q(x)=\frac{10\cos\left(  13x\right)  }{\left(  x+0.1\right)  ^{2}}+\pi
e^{x}\sin\left(  20.23x\right)  i,\quad x\in\left[  0,1\right]
\]
with the input data (\ref{inpdata}) of the following form. $101$ points
$\rho_{k}$ are chosen uniformly distributed on the segment $\left[
0.1,100\right]  $, and $a(\rho_{k})=\sin\rho_{k}$, $b(\rho_{k})=\cos\rho_{k}$.
The\ corresponding numbers $\ell_{k}$ were computed with the aid of the series
representations (\ref{phiNSBF}) and (\ref{S}) (for the details regarding the
computation of the coefficients $g_{n}(x)$ and $s_{n}(x)$ we refer to
\cite{KNT} and \cite[Sect. 9.4]{KrBook2020}). Now, having the input data we
solve the inverse Problem A following the algorithm from Section
\ref{Sect Summary}. In the first step, system (\ref{FirstSystem}) was solved
with $N=18$ (that gives $20$ computed coefficients for each of the
representations (\ref{phi1}) and (\ref{S1})). Here $\omega(L)$ was obtained
with the absolute error $2.49\cdot10^{-7}$, the absolute errors of the
computed $q(0)$ and $q(L)$ were $0.00084$ and $0.0012$, respectively.

Next, according to step 2 we chose a set of $700$ points $\left\{  \gamma
_{k}\right\}  _{k=1}^{700}$ logarithmically equally spaced on the segment
$\left[  0.1,1500\right]  $ and computed $\left\{  \varphi_{18}(\gamma
_{k},1),\,S_{18}(\gamma_{k},1)\right\}  _{k=1}^{700}$. The obtained values
were used to compute coefficients of the system (\ref{main system reduced}).
Solution of (\ref{main system reduced}) with the posterior recovery of $q(x)$
(subsection \ref{subsect recovery q}, option 1) gave us the result depicted in
Fig. 2. The absolute error of the computed potential was $3.5\cdot10^{-3}$
(option 1 from subsection \ref{subsect recovery q}) and $0.8\cdot10^{-3}$
(option 2). In both cases this is a remarkable accuracy for a rapidly
oscillating complex valued potential growing up to $10^{3}$.%

\begin{figure}
[ptb]
\begin{center}
\includegraphics[
height=3.6758in,
width=4.8941in
]%
{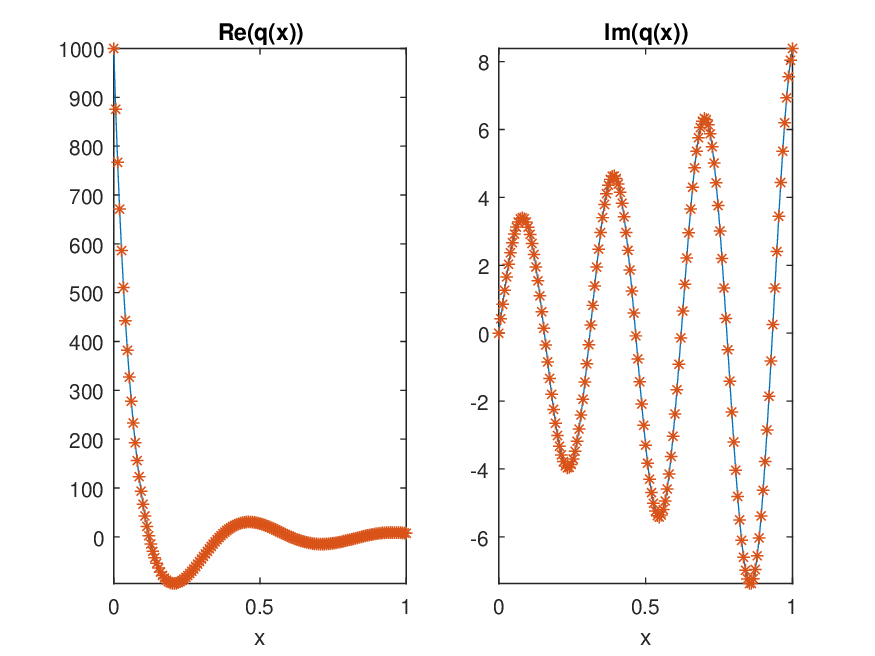}%
\caption{Real and imaginary parts of the recovered potential from Example 2.}%
\label{Fig2}%
\end{center}
\end{figure}

\textbf{Example 3. }Consider Problem A for the potential
\[
q(x)=\frac{1}{4}\left(  \left(  6x-\pi\right)  ^{6}-8\left(  6x-\pi\right)
^{4}+\left(  10.8x-\pi\right)  ^{2}\right)  +20.23+i\Gamma\left(
x+\pi\right)  ,\quad x\in\left[  0,1\right]  .
\]
With the same choice of the sets of points and parameters involved as in the
previous example, the result of the solution of the problem is shown in Fig.
3. The absolute error of the computed potential was $1.8\cdot10^{-4}$ (option
1 from subsection \ref{subsect recovery q}) and $2.4\cdot10^{-3}$ (option 2).
$\ $%

\begin{figure}
[ptb]
\begin{center}
\includegraphics[
height=3.4962in,
width=4.6552in
]%
{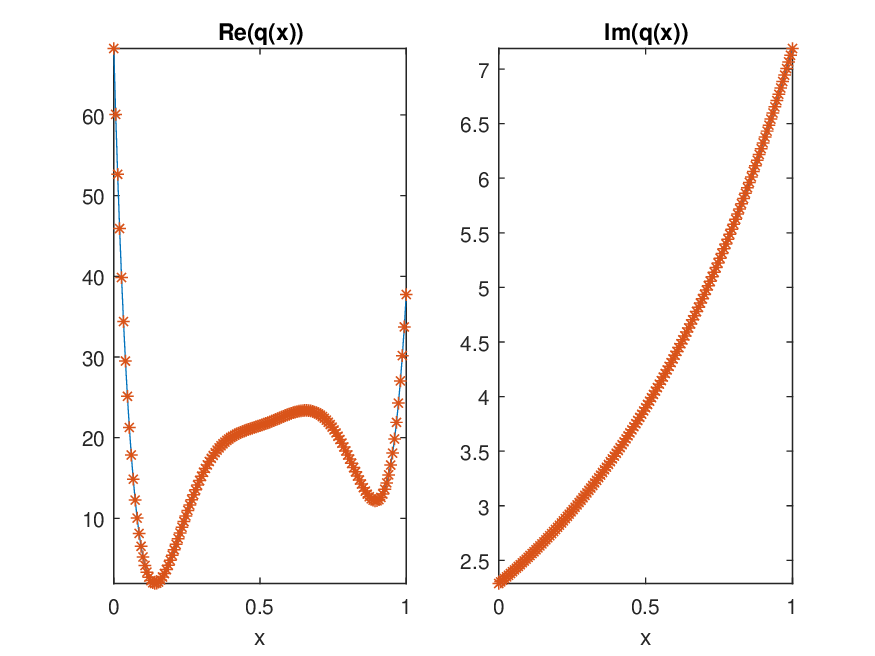}%
\caption{Real and imaginary parts of the recovered potential from Example 3.}%
\label{Fig3}%
\end{center}
\end{figure}

\textbf{Example 4. }Consider Problem A for the potential%
\[
q(x)=\int_{0}^{x}\left(  \left\vert s-\frac{1}{3}\right\vert +\pi\left\vert
s-\frac{4}{5}\right\vert \right)  ds+(1-\left(  \pi x-1\right)  ^{2}%
\operatorname*{sgn}(1-\pi x))i,\quad x\in\left[  0,1\right]  .
\]
Both the real and the imaginary parts of $q(x)$ are $C^{1}$-functions
possessing discontinuous second derivatives. The input data of the problem are
computed in $101$ uniformly distributed random points $\rho_{k}\in\left(
0,15\right)  $, with $a(\rho_{k})=\sin\rho_{k}$, $b(\rho_{k})=\cos\rho_{k}$.

In the first step, system (\ref{FirstSystem}) was solved with $N=18$ (which
means $20$ computed coefficients for each of the representations (\ref{phi1})
and (\ref{S1})). Here $\omega(L)$ was obtained with the absolute error
$5.33\cdot10^{-6}$, the absolute errors of the computed $q(0)$ and $q(L)$ were
$0.0007$ and $0.0014$, respectively. Next, according to step 2 we chose a set
of $700$ points $\left\{  \gamma_{k}\right\}  _{k=1}^{700}$ logarithmically
equally spaced on the segment $\left[  0.1,1500\right]  $ and computed
$\left\{  \varphi_{18}(\gamma_{k},1),\,S_{18}(\gamma_{k},1)\right\}
_{k=1}^{700}$. The obtained values were used to compute coefficients of the
system (\ref{main system reduced}). Solution of (\ref{main system reduced})
with the posterior recovery of $q(x)$ gave us the result depicted in Fig. 4.
The maximum absolute error of the potential (computed following both options
from subsection \ref{subsect recovery q}) resulted in $0.0025$.%

\begin{figure}
[ptb]
\begin{center}
\includegraphics[
height=3.1475in,
width=4.1949in
]%
{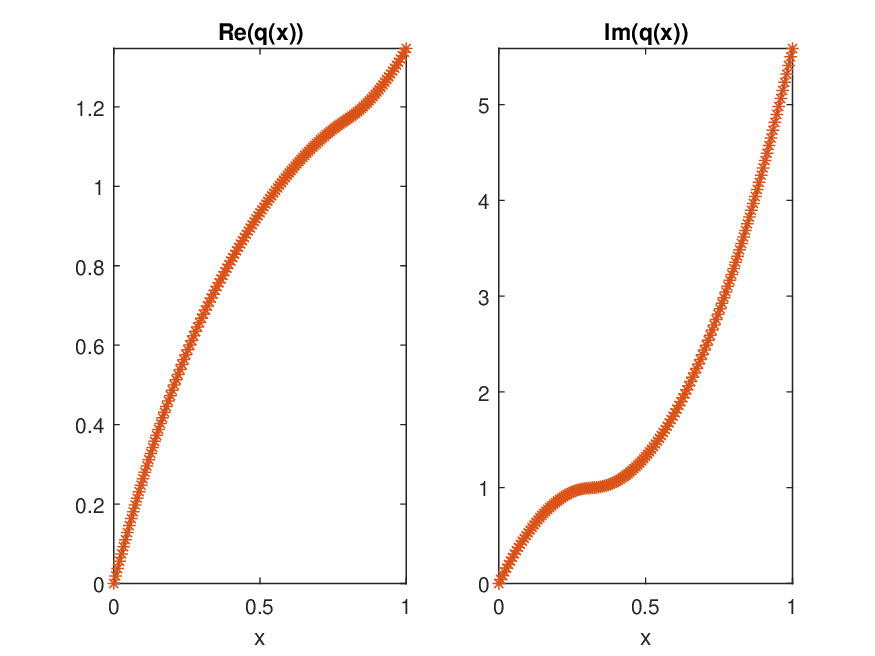}%
\caption{Real and imaginary parts of the recovered potential from Example 4.}%
\label{Fig4}%
\end{center}
\end{figure}

\section{Conclusions\label{Sect Conclusions}}

An approach for the numerical solution of a variety of inverse coefficient
problems for the Sturm-Liouville equation $-y^{\prime\prime}+q(x)y=\rho^{2}y$
with a complex valued potential $q(x)$ is presented. It reduces the solution
of a problem to the solution of two systems of linear algebraic equations for
the coefficients of the Neumann series of Bessel functions representations for
solutions of the Sturm-Liouville equation, and the potential $q(x)$ is
recovered from the first two components of the solution vector. The method is
easy for implementation, direct, accurate and applicable to a large variety of
inverse problems. Its performance is illustrated by numerical
examples.\bigskip

\textbf{Funding }Research was supported by CONAHCYT, Mexico via the project 284470.

\textbf{Data availability} The data that support the findings of this study
are available upon reasonable request.

\textbf{Conflict of interest }This work does not have any conflict of interest.

\end{document}